\pdfoutput=1
\RequirePackage{ifpdf}
\ifpdf 
\documentclass[pdftex]{sigma}
\else
\documentclass{sigma}
\fi

\usepackage{mathrsfs}

\numberwithin{equation}{section}

\newtheorem{prop}{Proposition}[section]
\newtheorem{theo}[prop]{Theorem}
\newtheorem*{theo*}{Theorem}

\newtheorem{conj}[prop]{Conjecture}

\theoremstyle{definition}

\newtheorem{quest}[prop]{Question}
\newtheorem{rema}[prop]{Remark}
\newtheorem{defi}[prop]{Definition}

\newcommand{\HH}{\mathbf{H}}

\newcommand{\RR}{\mathbf{R}}

\newcommand{\ZZ}{\mathbf{Z}}

\newcommand{\cH}{\mathcal H}

\DeclareMathOperator{\spt}{spt}

\DeclareMathOperator{\Ric}{Ric}

\DeclareMathOperator{\secondfund}{II}

\newcommand{\pa}[2]{\frac{\partial #1}{\partial #2}}

\newcommand{\bangle}[1]{\left\langle #1 \right\rangle}

\begin{document}
\allowdisplaybreaks

\newcommand{\arXivNumber}{2007.12563}

\renewcommand{\thefootnote}{}

\renewcommand{\PaperNumber}{099}

\FirstPageHeading

\ShortArticleName{Dihedral Rigidity of Parabolic Polyhedrons in Hyperbolic Spaces}

\ArticleName{Dihedral Rigidity of Parabolic Polyhedrons\\ in Hyperbolic Spaces\footnote{This paper is a~contribution to the Special Issue on Scalar and Ricci Curvature in honor of Misha Gromov on his 75th Birthday. The full collection is available at \href{https://www.emis.de/journals/SIGMA/Gromov.html}{https://www.emis.de/journals/SIGMA/Gromov.html}}}

\Author{Chao LI}

\AuthorNameForHeading{C.~Li}

\Address{Department of Mathematics, Princeton University,\\ Fine Hall, 304 Washington Rd, Princeton, NJ 08544, USA}
\Email{\href{chaoli@math.princeton.edu}{chaoli@math.princeton.edu}}
\URLaddress{\url{https://web.math.princeton.edu/~chaoli/}}

\ArticleDates{Received July 27, 2020, in final form September 30, 2020; Published online October 06, 2020}

\Abstract{In~this note, we establish the dihedral rigidity phenomenon for a collection of~parabolic polyhedrons enclosed by horospheres in hyperbolic manifolds, extending Gromov's comparison theory to metrics with negative scalar curvature lower bounds. Our~result is a localization of~the positive mass theorem for asymptotically hyperbolic manifolds. We~also motivate and formulate some open questions concerning related rigidity phenomenon and convergence of~metrics with scalar curvature lower bounds.}


\Classification{53C21; 53A10}

\begin{flushright}
\begin{minipage}{80mm}
\it Dedicate this paper to Professor Misha Gromov \\
on the occasion of his 75th birthday.
\end{minipage}
\end{flushright}

\renewcommand{\thefootnote}{\arabic{footnote}}
\setcounter{footnote}{0}

\section{Introduction}
In~\cite{Gromov2014Dirac}, Gromov proposed a geometric comparison theory for metrics with scalar curvature lower bounds. He speculated that Riemannian polyhedrons should play the role of~triangles in Alexandrov's comparison theory for sectional curvature~\cite{AleksandrovRerestovskiiNikolaev86}. As a first step, he obtained the following theorem for metrics with nonnegative scalar curvature, where the comparison models are Euclidean cubes:

\begin{theo}[\cite{Gromov2014Dirac}]\label{theo.cube.nonrigit}
Let $M=[0,1]^n$ be a cube, $g$ a smooth Riemannian metric. Then $(M,g)$ cannot simultaneously satisfy:
\begin{enumerate}\itemsep=0pt
\item[${\rm 1)}$] the scalar curvature $R(g)\ge 0$;
\item[${\rm 2)}$] each face of~$M$ is weakly strictly mean convex;\footnote{In~this paper, the mean curvature is taken with respect to outer unit normal vector. For instance, the standard sphere $S^{n-1}$ in $\RR^n$ has mean curvature $n-1$.}
\item[${\rm 3)}$] the dihedral angles between adjacent faces are all acute.
\end{enumerate}
\end{theo}

Theorem~\ref{theo.cube.nonrigit} also has a rigidity statement: if $n\le 7$, and we assume all dihedral angles are not larger than $\pi/2$ in condition (3), then $(M,g)$ is isometric to an Euclidean rectangular solid (see~\cite{Li2017polyhedron,li2019dihedral}). This is called the dihedral rigidity phenomenon. In~\cite{Gromov2014Dirac,Gromov2018Adozen}, Gromov conjectured that this property is satisfied for all convex polyhedrons in $\RR^n$:

\begin{conj}[the dihedral rigidity conjecture]\label{conj.dihedral.rigidity}
Let $M\subset \RR^n$ be a convex polyhedron and~$g_0$ be the Euclidean metric. Suppose $g$ is a smooth Riemannian metric on~$M$. Denote its faces by~$F_i$, the mean curvature of~$F_i$ by~$H_i$, and the dihedral angle between two adjacent faces $F_i$, $F_j$ by~$\measuredangle_{ij}$ $(\measuredangle_{ij}(g)$ may be nonconstant$)$. Assume:
\begin{enumerate}\itemsep=0pt
\item[${\rm 1)}$] $R(g)\ge 0$ in $M$;
\item[${\rm 2)}$] $H_i(g)\ge 0$ on each face $F_i$;
 \item[${\rm 3)}$] $\measuredangle_{ij}(g)\le \measuredangle_{ij}(g_0)$ on each pair of~adjacent faces $F_i$, $F_j$.
\end{enumerate}
Then $(M,g)$ is isometric to a flat polyhedron in $\RR^n$.
\end{conj}

Conjecture~\ref{conj.dihedral.rigidity} and related problems have been studied and extended in recent years (see, e.g.,~\cite{Gromov2018Metric,gromov2018IAS,gromov2019lectures,Li2017polyhedron,li2019dihedral,LiMantoulidis2018positive,miao2019measuring}), leading to a range of~interesting new discoveries and questions on~manifolds with nonnegative scalar curvature. In~this paper, we investigate the analogous polyhedral comparison principle, together with the rigidity phenomenon, for metrics with nega\-tive scalar curvature lower bound.

By scaling, we assume $R(g)\ge -n(n-1)$. Our comparison model is a collection of~polyhedrons in the hyperbolic space, called \textit{parabolic prisms}, which we define now. Let $(\HH^n,g_H)$ be the hyperbolic space with sectional curvature $-1$. We~choose the coordinate system $\{x_1,\dots,x_n\}$, $x_j\in \RR$, such that $g_H$ takes the form
\[g_H={\rm d}x_1^2+{\rm e}^{2x_1}\big({\rm d}x_2^2+\cdots+{\rm d}x_n^2\big).\]
For any constant~$c$, the coordinate hyperplane $x_1=c$ is umbilical with constant mean curvature $n-1$ with respect to~$\partial_{x_1}$. The induced metric on it is isometrically Euclidean. These hyperplanes are called \textit{horospheres}. For $j\ge 2$, the coordinate hyperplanes $x_j=c$ are totally geodesic, and they intersect each other and the horospheres orthogonally.

Denote $\hat{x}=(x_2,\dots,x_n)$. Given a polyhedron~$P\subset \RR^{n-1}$, we call the set $\{(x_1,\hat{x})\colon 0\le x_1\le 1,\, \hat{x}\in P\}$ a parabolic prism in $\HH^n$. As a special case of~the main theorem of~this paper, Theorem~\ref{theo.main}, we have the following
\begin{theo}\label{theo.parabolic.cube}
Let $n\le 7$, $M=[0,1]^n$ be a parabolic rectangle in $\HH^n$, $g_H$ be the hyperbolic metric on~$M$. Denote the face $\partial M\cap \{x_1=1\}$ by~$F_T$, the face $\partial M\cap \{x_1=0\}$ by~$F_B$. Assume~$g$ is a Riemannian metric on~$M$ such that:
\begin{enumerate}\itemsep=0pt
\item[${\rm 1)}$] $R(g)\ge -n(n-1)$ in $M$;
\item[${\rm 2)}$] $H(g)\ge n-1$ on~$F_T$, $H(g)\ge -(n-1)$ on~$F_B$, and~$H(g)\ge 0$ on~$\partial M\setminus (F_T\cup F_B)$;
\item[${\rm 3)}$] the dihedral angles between adjacent faces of~$M$ are everywhere not larger than~$\pi/2$.
\end{enumerate}
Then $(M,g)$ is isometric to a parabolic rectangle in $\HH^n$.
\end{theo}

Theorem~\ref{theo.parabolic.cube} addresses a question discussed during the workshop ``Emerging Topics on Scalar Curvature and Convergence'' at IAS in October, 2018. See \cite[Section~6]{gromov2018IAS}.

We~remark that Theorem~\ref{theo.parabolic.cube} holds for any general Riemannian polyhedron~$(M,g)$ with a~proper polyhedral map to the parabolic cube of~nonzero degree. It also holds for more general polyhedral types, as long as the comparison model is a parabolic prism $[0,1]\times P$, and~$P\subset \RR^{n-1}$ satisfies Conjecture~\ref{conj.dihedral.rigidity}. By~\cite{Li2017polyhedron,li2019dihedral}, $P$ can be any $3$-dimension simplices or $n$-dimensional non-obtuse prisms. See Theorem~\ref{theo.dihedral.rigidity}.

It has been observed in \cite[Section~5]{li2019dihedral} that Conjecture~\ref{conj.dihedral.rigidity} is a localization of~the positive mass theorem for asymptotically flat manifolds \cite{SchoenYau1979ProofPositiveMass,Witten1981newproof}. Analogously, Theorem~\ref{theo.parabolic.cube} localizes the positive mass theorem for asymptotically hyperbolic manifolds (see, e.g.,~\cite{chruscielHerzlich2003mass,ChruscielNagy2001mass,Wang2001mass}), which, in special cases, can be deduced from the following rigidity result of~scalar curvature, due to Min-Oo~\cite{Minoo1989rigidity} on spin manifolds and to Andersson--Cai--Galloway~\cite{AnderssonCaiGalloway2008rigidity} on all manifolds of~dimension at most $7$ (see also~\cite{chruciel2019hyperbolic,HuangJangMartin2020mass} for more recent developments):
\begin{theo}[\cite{AnderssonCaiGalloway2008rigidity}; see also~\cite{Minoo1989rigidity} for spin manifolds]\label{theo.R.rigidity.hyperbolic}
Suppose $(S^n,g)$, $2\le n\le 7$, has scalar curvature $R(g)\ge -n(n-1)$, and is isometric to $\HH^n$ outside a compact set. Then $(S^n,g)$ is isometric to $(\HH^n, g_H)$.
\end{theo}

Precisely, suppose Theorem~\ref{theo.parabolic.cube} holds. Given $(S^n,g)$ with $R(g)\ge -n(n-1)$ such that $(S,g)$ is isometric to $(\HH^n,g_H)$ outside a compact set $K$, take a sufficiently large $R$ such that the boundary of~the parabolic rectangle $M=[-R,R]^n$ isometrically embeds into $S\setminus K$. Denote the region bounded by~$\partial M$ in $S$ by~$M_1$. Then $M_1$ has a degree one map to $M$, by sending $K$ to an interior point $p\in M$ and~$M_1\setminus K$ to $M\setminus \{p\}$. Thus Theorem~\ref{theo.parabolic.cube} implies that $(M_1,g)$ is isometric to $(M,g_H)$.

Given the connection between the positive mass theorem and the dihedral rigidity conjecture, it would be interesting to see whether a similar comparison principle holds for metrics with positive scalar curvature lower bound. The delicate issue is that the corresponding rigidity phenomenon on a hemisphere is false, due to the counterexamples by Brendle--Marques--Neves~\cite{Brendle2011deformation}. This, together with other related open questions, will be discussed in Section~\ref{section4}.

\section{Notations and the main theorem}
The main objects in this paper are \textit{Riemannian polyhedrons}, which we define as follows.

\begin{defi}A compact Riemannian manifold $(M^n,g)$ with boundary is called a Riemannian polyhedron, if $(M,g)$ can be isometrically embedded into $\RR^N$ for some $N\ge n$, and at every $x\in M$, there is a radius $r>0$ and a diffeomorphism $\phi_x\colon B_r\big(x\in \RR^N\big)\rightarrow B_1\big(0^N\big)$, such that $\phi_x(B_r\cap M)=P\cap B_1\big(0^N\big)$ for some Euclidean polyhedral cone~$P$ of~dimension~$n$, and~$D\phi_x|_x$ is an isometry. Further, we require that $\phi_x$ is $C^{2,\alpha}$ for some $\alpha\in (0,1)$ independent of~$x$.
\end{defi}

Specially, a compact domain $M$ enclosed by piecewise $C^{2,\alpha}$ hypersurfaces in a smooth Riemannian manifold is a Riemannian polyhedron. Given $x\in M^n$, there is an integer $k\in [0, n]$ such that a neighborhood of~$x$ in $M$ is diffeomorphic to $P_0^{n-k}\times \RR^k$, and~$P_0$ is a polyhedral cone in $\RR^{n-k}$ without translation symmetry. We~call the union of~all such points the $k$-faces of~$M$. In~particular, the $n$-face is the interior of~$M$, the $(n-1)$-faces are the union of~smooth components of~$\partial M$ (which we called ``faces'' in Theorem~\ref{theo.parabolic.cube}), and the $(n-2)$-faces are the interior of~edges of~$M$.

\begin{defi}Let $P\subset \RR^n$ be a flat Euclidean polyhedron, and~$(M^n,g)$ be a Riemannian poly\-hedron. We~say $M$ is over-$P$-polyhedral, if $M$ admits a proper polyhedral map $\phi$ onto $P$ (i.e., $\phi$ maps any $k$-face of~$M$ to a $k$-face of~$P$), such that $\phi$ is of~nonzero degree.
\end{defi}

In~\cite{Li2017polyhedron} and~\cite{li2019dihedral}, Conjecture~\ref{conj.dihedral.rigidity} was proved for a Riemannian polyhedrons that are over-$P$-polyhedral, where
\begin{enumerate}\itemsep=0pt
\item[1)] either $n=3$, and~$P\subset\RR^3$ is an arbitrary simplex;
\item[2)] or $3\le n\le 7$, and~$P$ is the Cartesian product $P_0^2\times [0,1]^{n-2}$. Here $P_0\subset \RR^2$ is a polygon with non-obtuse dihedral angles.
\end{enumerate}
Precisely, we have:

\begin{theo}[\cite{Li2017polyhedron,li2019dihedral}]\label{theo.dihedral.rigidity}
Let $P\subset \RR^n$ be as above, $(M^n,g)$ be an over-$P$-polyhedral Riemannian polyhedron, and~$\phi\colon M\rightarrow P$ be the polyhedral map of~nonzero degree. Suppose:
\begin{enumerate}\itemsep=0pt
\item[${\rm 1)}$] $R(g)\ge 0$ in $M$;
\item[${\rm 2)}$] $H(g)\ge 0$ on each face of~$M$;
\item[${\rm 3)}$] $\measuredangle_{ij}(g)|_x\le \measuredangle_{ij}(g_0)|_{\phi(x)}$ for every point~$x$ on the edges of~$M$.
\end{enumerate}
Then $(M,g)$ is isometric to an Euclidean polyhedron.
\end{theo}

We~now state the main theorem of~this paper.

\begin{theo}\label{theo.main}Let $2\le n \le 7$, $P\subset \RR^{n-1}$ be an Euclidean polyhedron such that Theorem~{\rm \ref{theo.dihedral.rigidity}} holds for $P$. Denote $([0,1]\times P,g_H)$ the parabolic prism in the hyperbolic space with $R(g)=-n(n-1)$. Suppose $(M^n,g)$ is a Riemannian polyhedron that is over-$[0,1]\times P$-polyhedral, and~$\phi\colon M^n\rightarrow [0,1]\times P$ be the polyhedral map of~nonzero degree. Suppose:
\begin{enumerate}\itemsep=0pt
\item[${\rm 1)}$] $R(g)\ge -n(n-1)$ in $M$;
\item[${\rm 2)}$] $H(g)\ge n-1$ on~$\phi^{-1}(\{1\}\times P)$, $H(g)\ge -(n-1)$ on~$\phi^{-1}(\{0\}\times P)$, and~$H(g)\ge 0$ on other faces of~$M$;
\item[${\rm 3)}$] $\measuredangle_{ij}(g)|_x\le \measuredangle_{ij}(g_H)|_{\phi(x)}$ for every point $x$ on the edges of~$M$.
\end{enumerate}
Then $(M,g)$ is isometric to a parabolic prism in the hyperbolic space.
\end{theo}

\begin{rema}
The dimension restriction~$n\le 7$ in Theorems~\ref{theo.dihedral.rigidity} and~\ref{theo.main} is due to regularity of~free boundary area minimizing surfaces and isoperimetric regions. In~light of~the recent progress on positive mass theorem in higher dimensions~\cite{SchoenYau2017positive}, we speculate the singular analysis may be applicable a non-rigid variance of~Theorem~\ref{theo.main}, i.e.,~Theorem~\ref{theo.cube.nonrigit}.
\end{rema}

\begin{rema}
The condition that $P=P_0\times \RR^{n-2}$, $P_0$ is non-obtuse, is due to the boundary regularity theory developed by Edelen and the author~\cite{EdelenLi2020regularity}. In~general, one can guarantee that free boundary area minimizing surfaces are $C^{2,\alpha}$ regular in a non-obtuse polyhedral domain. See~\cite[Section~9]{EdelenLi2020regularity}.
\end{rema}

\section{Proof of~the main theorem}

The proof of~Theorem~\ref{theo.main} is an adaptation of~the proof of~Theorem~\ref{theo.dihedral.rigidity}. We~will be using lots of~techniques developed in~\cite{li2019dihedral}. Given a Riemannian polyhedron~$(M^n,g)$ as in Theorem~\ref{theo.main}, denote the faces $F_T=\phi^{-1}(\{1\}\times P)$ and~$F_B=\phi^{-1}(\{0\}\times P)$. If there exists a point $x$ on the edge of~$M$ such that $\measuredangle(g)|_x<\measuredangle(g_H)|_{\phi(x)}$, we deform the metric $g$ to $\tilde{g}$ as in \cite[Section~11]{Gromov2018Metric}, such that $\measuredangle(\tilde{g})|_x=\measuredangle(g_H)|_{\phi(x)}$ in a neighborhood of~$x$ and the mean curvature of~the two faces containing $x$ increases. Thus, without loss of~generality, we assume that $\measuredangle(g)|_x=\measuredangle(g_H)|_{\phi(x)}$ for~all points $x$ on the edge.

Consider the relative isoperimetric problem:
\begin{gather}
I=\inf\big\{\cH^{n-1}\big(\partial \Omega\llcorner \mathring{M}\big)-(n-1)\cH^n(\Omega)\colon
\Omega\subset M \text{ is a Caccioppoli set},\nonumber
\\ \phantom{I=\inf\{}
F_B\subset \Omega,\, F_T\cap \Omega=\varnothing\big\}.\label{problem.variation}
\end{gather}

It follows from the standard compactness results that $I$ is achieved by a Caccioppoli set $\Omega$. Denote $\Sigma=\spt\big(\partial \Omega\llcorner \mathring{M}\big)$. Since $F_B$, $F_T$ meet other faces orthogonally and~$H(g)\ge n-1$ on~$F_T$, $H(g)\ge -(n-1)$ on~$F_B$, by the strong maximum principle (see \cite[Section~3.1]{li2019dihedral}), either $\Sigma$ is disjoint from $F_T$ and~$F_B$, or $\Sigma$ coincides with $F_T$ or $F_B$. In~any case, $\Sigma$ is an isoperimetric surface with free boundary on~$\partial M\setminus (F_T\cup F_B)$.

We~remark here that similar variational problems as \eqref{problem.variation} have been considered by Witten--Yau~\cite{WittenYau1999connectedness} and by Andersson--Cai--Galloway~\cite{AnderssonCaiGalloway2008rigidity}, where it is called the BPS brane action.

We~now study the regularity of~$\Sigma$. Since $n\le 7$, the regularity of~$\Sigma$ in the interior and smooth part of~$\partial M$ follows from the classical theory~\cite{GruterJost1986Allard, Simons1968minimal}. For a point $x\in \Sigma$ and a $k$-face of~$M$ with $k\le n-2$, we note that the tangent domain of~$M$ at $x$ is given by~$W^2\times [0,\infty)^{k-2}\times \RR^{n-k}$, where~$W^2$ is a wedge region in $\RR^2$ with non-obtuse opening angle. Thus, we apply \cite[Section~9]{EdelenLi2020regularity} and \cite[Appendix B]{li2019dihedral}, and conclude:
\begin{prop}
$\Sigma$ is $C^{2,\alpha}$ graphical over its tangent plane everywhere.
\end{prop}

Moreover, since $\Sigma$ is homologous to $F_B$, we conclude that at least one connected component (which we still denote by~$\Sigma$) has a nonzero degree map to $P$, given by~$\Sigma\simeq F_B\xrightarrow{\phi} \{0\}\times P$.

Since $\Omega$ is a minimizer for \eqref{problem.variation}, $\Sigma$ has constant mean curvature $(n-1)$ with respect to the outer unit normal $\nu$ of~$\Omega$, and stability implies that
\begin{gather}
Q(\varphi):=\int_\Sigma |\nabla \varphi|^2 -\frac{1}{2}\big(R_M-R_\Sigma + n(n-1)+|\mathring{A}|^2\big)\varphi^2 {\rm d}\cH^{n-1}\nonumber
\\ \hphantom{Q(\varphi):=}
{} -\int_{\partial \Sigma} \secondfund(\nu,\nu)\varphi^2 {\rm d}\cH^{n-2}\ge 0,
\end{gather}
for all $C^1$ function~$\varphi$. Here $R_\Sigma$, $\mathring{A}$ are the scalar curvature of~the induced metric and the traceless second fundamental of~$\Sigma$, respectively, and~$\secondfund$ is the second fundamental form of~$\partial M$.

Let $\varphi>0$ be the principal eigenfunction associated with the quadratic form $Q$. Then $\varphi$ solves the equation
\begin{equation}
\begin{cases}
\Delta_\Sigma \varphi+\dfrac{1}{2}\big(R_M-R_\Sigma + n(n-1)+|\mathring{A}|^2\big)\varphi = -\lambda_1\varphi\,,
\\
\displaystyle\pa{\varphi}{\eta}=\secondfund(\nu,\nu)\varphi\,.
\end{cases}
\end{equation}
Here $\eta$ is the outer conormal vector field of~$\Sigma$, and~$\lambda_1$ is the principal eigenvalue associated with~$Q$. It follows from \cite[Lemma 4.1]{li2019dihedral} that $\varphi\in C^{2,\alpha}(\Sigma)$. Denote $\tilde{g}=\varphi^{\frac{2}{n-2}}g$ on~$\Sigma$. By the very same calculations as in \cite[equations~(4.6) and~(4.7)]{li2019dihedral}, we have
\[R(\tilde{g})=\varphi^{-\frac{n}{n-2}}\left(\big(R_M+n(n-1)+|\mathring{A}|^2+\lambda_1\big)\varphi+\frac{n-1}{n-2}\frac{|\nabla \varphi|^2}{\varphi}\right)\ge 0,\]
and~$H_{\partial \Sigma}(\tilde g)=\varphi^{-\frac{1}{n-2}}(H_{\partial \Sigma}(g)+\secondfund(\nu,\nu)=\varphi^{-\frac{1}{n-2}} H_{\partial M}(g)\ge 0$.

Moreover, since $\Sigma$ meets $\partial M$ orthogonally and~$\tilde{g}$ is conformal to $g$, the dihedral angles of~$(\Sigma,\tilde{g})$ is everywhere equal to that of~$P$. Thus, by Theorem~\ref{theo.dihedral.rigidity}, $(P,\tilde{g})$ is isometric to an Euclidean polyhedron. Tracing equality, we have
\[
R_M=0,\qquad \mathring A=0,\qquad \lambda_1=0,\qquad \nabla \varphi=0 \qquad \text{on }\Sigma.
\]
Therefore $\varphi$ is a constant function, and hence $\Ric_M(\nu,\nu)=-(n-1)$ on~$\Sigma$ and~$\secondfund(\nu,\nu)=0$ on~$\partial \Sigma$. It follows that $\Sigma$ is totally umbilical and infinitesimally rigid.

Next, we adapt the ideas in~\cite{CarlottoChodoshEichmair2016effective,ChodoshEichmairMoraru2018splitting} to study rigidity. Let $M^-$ be the region enclosed by~$\Sigma$ and~$F_B$. We~follow the same argument as in \cite[Section~4]{li2019dihedral}: by constructing the very same deformed metrics $\{g(t)\}_{t\in [0,\varepsilon)}$, solving the relative isoperimetric problem \eqref{problem.variation} and taking convergence as $t\rightarrow 0$, we obtain another free-boundary isoperimetric hypersurface $\Sigma'$ in $(M,g)$ lying between $\Sigma$ and~$F_B$. Moreover, $\Sigma'$ is also isometrically Euclidean, and is umbilical and infinitesimally rigid. By repeating this argument, we obtain a dense collection of~such hypersurfaces~$\{\Sigma^\rho\}$ in~$M$.

Fix $\Sigma^\rho$, its outer unit normal $\nu$, and~$x_0\in \Sigma^\rho$. For $\rho_j$ sufficiently close to $\rho$, $\Sigma^{\rho_j}$ can be writ\-ten as a normal graph of~function~$u^j$ over $\Sigma^\rho$. By standard curvature, the function~$u^j/u^j(x_0)$ converges in $C^{2,\alpha}(\Sigma^\rho)$ to a nonzero function~$u$. The Gauss--Codazzi equation implies that, for~any tangential vector $X,Y$ on~$\Sigma^\rho$,
\[
\big(\nabla_{\Sigma^\rho}^2 u\big)(X,Y)+Rm_M(\nu,X,Y,\nu)u + A_{\Sigma^\rho}(X,Y)u=0.
\]
Taking trace, we have that $\Delta_{\Sigma^\rho}u=0$. Also, since $\Sigma^{\rho_j}$ meets $\partial M$ orthogonally and~$\secondfund(\nu,\nu)=0$ on~$\partial\Sigma^{\rho_j}$, $\pa{u}{\eta}=0$. Thus $u$ is a constant function, and hence $Rm_M(\nu,X,Y,\nu)=-\bangle{X,Y}$. This proves that $M$ has constant sectional curvature $-1$. Theorem~\ref{theo.main} is proved.

\section{Discussions and related questions}\label{section4}
\subsection{Metrics with positive scalar curvature lower bounds}
It is tempting to conjecture that a suitable extension of~Theorems~\ref{theo.dihedral.rigidity} and~\ref{theo.main} holds for metrics with positive scalar curvature lower bounds, where the model space is $S^n$ with a~round metric. Although the precise formulation is unclear, such an extension will likely to localize certain scalar curvature rigidity phenomenon for spheres. Recall the following theorem by Brendle and Marques~\cite{BrendleMarques2011scalar}:

\begin{theo}[\cite{BrendleMarques2011scalar}]
Let $\Omega=B(\delta)\subset S^n$ be a closed geodesic ball of~radius $\delta$ with $\cos\delta\ge \frac{2}{\sqrt{n+3}}$, and~$\overline{g}$ be the standard metric on~$S^n$. Suppose $g$ is another metric on~$\Omega$ such that:
\begin{enumerate}\itemsep=0pt
\item[${\rm 1)}$] $R(g)\ge R(\overline{g})$ in $\Omega$;
\item[${\rm 2)}$] $H(g)\ge H(\overline{g})$ on~$\partial \Omega$;
\item[${\rm 3)}$] $g$ and~$\overline{g}$ induce the same metric on~$\partial \Omega$.
\end{enumerate}
If $g-\overline{g}$ is sufficiently small in the $C^2$-norm, then $g$ is isometric to $\overline{g}$.
\end{theo}

The lower bound $\frac{3}{\sqrt{n+3}}$ for $\delta$ was improved in a subsequent paper by Cox, Miao and~Tam~\cite{CoxMiaoTam2013remarks}. However, it is known that the analogous statement for $\delta=\frac{\pi}{2}$ does not hold, due to the counterexample by Brendle, Marques and Neves~\cite{Brendle2011deformation}.

The original proof of~Theorem~\ref{theo.cube.nonrigit} by Gromov uses the fact that a cube is the fundamental domain of~$\ZZ^n$ action on~$\RR^n$: assuming a counterexample for Theorem~\ref{theo.cube.nonrigit} exists, through a~sequence of~doubling and smoothing, one obtains a smooth metric on~$T^n$ with positive scalar curvature, contradicting~\cite{GromovLawson1980Spin} and~\cite{SchoenYau1979structure}.

Take the standard embedding $S^n_+\subset\RR^{n+1}$. The hemisphere $S^n_+=S^n\cap \{x_{n+1}\ge 0\}$ can be obtained by consecutive doublings of~the spherical simplex\vspace{-.5ex}
\[
\Omega_n:=S^n\cap \{x_j\ge 0,\, j=1,\dots,n+1\}.
\]

We~make the following conjecture concerning dihedral rigidity of~$\Omega_n$.

\begin{conj}\label{conj.spherical.simplex}
Let $(M^n,g)$ be a Riemannian polyhedron which is diffeomorphic to a simplex of~dimension~$n$. Suppose\vspace{-.5ex}
\begin{enumerate}\itemsep=0pt
\item[${\rm 1)}$] $R(g)\ge n(n-1)$ in $M$;
\item[${\rm 2)}$] $\partial M$ is piecewise totally geodesic;
\item[${\rm 3)}$] $\measuredangle_{ij}(g)\le \frac{\pi}{2}$ on the edges of~$M$;
\item[${\rm 4)}$] moreover, each face of~$M$ is globally isometric to the standard spherical simplex $\Omega_{n-1}$.
\end{enumerate}
Then $(M^n,g)$ is isometric to $\Omega_n$.
\end{conj}

When $n=2$, Conjecture~\ref{conj.spherical.simplex} holds by doubling $M$ twice across its boundary, and using a~theorem due to Toponogov~\cite{Toponogov1959evaluation}. On the other hand, the construction in~\cite{Brendle2011deformation} does not seem to give a counterexample to Conjecture~\ref{conj.spherical.simplex}.

\subsection[Weak notions of $R \ge \kappa$]{Weak notions of~$\boldsymbol{R\ge \kappa}$}

One of~Gromov's motivations of~studying Conjecture~\ref{conj.dihedral.rigidity} is to define the notion of~``$R\ge 0$'' in~the weak sense. The crucial observation is that, the conditions (2), (3) concern $C^0$ properties of~the metric $g$, and are stable under $C^0$ convergence of~metrics (see~\cite{Gromov2014Dirac}). Thus, we may define{\samepage
\begin{gather*}
\text{``}R(g)\ge 0\text{''} \Leftrightarrow \text{there exists no cube }M \nonumber
\\ \hphantom{\text{``}R(g)\ge 0\text{''} \Leftrightarrow{}}
\text{with mean convex faces and everywhere acute dihedral angle.}
\end{gather*}}

And more generally for $\kappa<0$,
\begin{gather}
\text{``}R(g)\ge \kappa \text{''} \Leftrightarrow \text{there exists no cube }M \text{ in the product with }S^2_{-\kappa}\nonumber
\\ \hphantom{\text{``}R(g)\ge \kappa \text{''} \Leftrightarrow{}}
\text{with mean convex faces and everywhere acute dihedral angle.}\label{weak.notion.general}
\end{gather}
Here $S^2_{-\kappa}$ is the space form with scalar curvature $-\kappa$.

Using this observation, Gromov proved the following theorem on the convergence of~metrics with scalar curvature lower bounds:

\begin{theo}[\cite{Gromov2014Dirac}, see also~\cite{bamler2016ricci}]
Let $M^n$ be a smooth manifold, $g$, $g_k$, $k\ge 1$, be smooth Riemannian metrics on~$M$, and~$g_k\rightarrow g$ in $C^0$ as tensors. Suppose $R(g_k)\ge \kappa$ on~$M$. Then $R(g)\ge \kappa$ as well.
\end{theo}

Based on Theorem~\ref{theo.main}, we may define $R\ge \kappa$ for a negative constant $\kappa$:
\begin{gather}
\text{``}R(g)\ge -\kappa\text{''} \Leftrightarrow \text{there exists no cube }M \text{ with acute dihedral angles}\nonumber
\\ \hphantom{\text{``}R(g)\ge -\kappa\text{''} \Leftrightarrow{}}
\text{and faces }\{F_j\}, \text{such that } H>-(n-1) \text{ on }F_1,\nonumber
\\ \hphantom{\text{``}R(g)\ge -\kappa\text{''} \Leftrightarrow{}}
H>0 \text{ on other faces, and }H>(n-1) \text{ on the opposite face of~}F_1.\label{weak.notion.R>-kappa}
\end{gather}
\eqref{weak.notion.general} and \eqref{weak.notion.R>-kappa} should be equivalent on smooth metrics, but we~think that \eqref{weak.notion.R>-kappa} is sightly more natural conceptually, as it also satisfies the dihedral rigidity phenomenon.

Recently, Burkhardt-Guim~\cite{burkhardt-guim2019pointwise} proposed a different possible notion of~``$R>\kappa$'' using Ricci flow. See \cite[Definition~1.2]{burkhardt-guim2019pointwise}. These definitions all share some good properties as a weak notion. For instance, they all agree to $R(g)>\kappa$ in the classical sense for a $C^2$ metric $g$, and they can be localized in a neighborhood of~any point on the manifold. The natural question is:

\begin{quest}\label{question.different.definitions}
Given a smooth manifold $M$ and a $C^0$ metric $g$ on it. Do the definitions \eqref{weak.notion.R>-kappa} and \cite[Definition~1.2]{burkhardt-guim2019pointwise} agree on~$g$?
\end{quest}

\subsection*{Acknowledgements}

I would like to thank Christina Sormani and Misha Gromov for organizing the excellent workshop ``Emerging Topics on Scalar Curvature and Convergence'' at the Institute for Advanced Study, and everyone in the workshop for valuable discussions. The author is supported by NSF grant DMS-2005287.

\pdfbookmark[1]{References}{ref}
\LastPageEnding

\end{document}